\documentclass[twoside,reqno]{amsart} 

\usepackage{amsmath,amsthm,amssymb,amsxtra} 
\usepackage{latexsym,cite,eucal,rotate,url} 
\usepackage{float,multirow,enumitem} 
\usepackage{color,graphics,graphicx}

\newcommand{\fns}{\footnotesize}

\newcommand{\ang}[1]{\langle{#1}\rangle}

\newcommand{\mb}[1]{\mathbb{#1}}

\newcommand{\sbs}[1]{_{\substack{#1}}}

\newcommand{\hyq}[4]{\left[\setlength
\arraycolsep{.5mm}\ba{#1}#3\\#4\ea{\Big|#2}\right]}

\newcommand{\binm}{\binom}
\newcommand{\binq}[2]{\bigg[\genfrac{}{}{0mm}{0}{#1}{#2}\bigg]}

\newcommand{\mult}[2]{\setlength
\arraycolsep{.5mm}\begin{array}{#1}#2\end{array}}%

\newcommand{\be}{\begin{equation}}
\newcommand{\ee}{\end{equation}}
\newcommand{\ba}{\begin{array}}
\newcommand{\ea}{\end{array}}
\newcommand{\bmn}{\begin{eqnarray}}
\newcommand{\emn}{\end{eqnarray}}
\newcommand{\bnm}{\begin{eqnarray*}}
\newcommand{\enm}{\end{eqnarray*}}
\newcommand{\bln}{\begin{subequations}}
\newcommand{\eln}{\end{subequations}}

\newcommand{\pq}[1]{\begin{equation}#1\end{equation}}
\newcommand{\pp}[2]{\begin{aligned}[#1]#2 
            \end{aligned}}  
\newcommand{\pnq}[1]{\begin{align*}#1
            \end{align*}}    

\newcommand{\centro}[1]
           {\begin{center}#1\end{center}}

\newcommand{\lam}{\lambda}

\newcommand{\sig}{\sigma}

\newcommand{\ii}{\imath}
\newcommand{\jj}{\jmath}

\newcommand{\Lam}{\Lambda}
\newcommand{\Ome}{\Omega}

\newtheorem{thm}{Theorem}
\newtheorem{lemm}[thm]{Lemma}
\newtheorem{corl}[thm]{Corollary}

\newcommand{\referxy}[4]{\bibitem{kn:#1}{#2,}~\emph{#3,}~{#4.}}	
\newcommand{\cito}[1]{\cite{kn:#1}}	
\newcommand{\citu}[2]{\cite[#2]{kn:#1}}

\begin{document}{\fbox{\fns\today}\hfill}
\title{Summation Formulae for Binomial Moments}
\author{Marta Na Chen and Wenchang Chu}
\address{School of Mathematics and Statistics\newline
         Zhoukou Normal University\newline
         Zhoukou (Henan), China}
\email{chennaml@outlook.com}
\address{Via Dalmazio Birago 9/E, Lecce 73100, Italy}
\email{hypergeometricx@outlook.com}
\thanks{Corresponding author (W.~Chu): hypergeometricx@outlook.com}


\begin{abstract}
By combining the telescoping method with an algebraic relation,
four classes of binomial moments are examined. Several explicit
summation formulae are established. 
%
%
\end{abstract}

\maketitle\thispagestyle{empty}\vspace*{-5mm} 
\markboth{Summation Formulae for Binomial Moments}
	 {M.~N.~Chen and W.~Chu}

\section{Introduction and Outline}

Let $\mb{N}$ stand for the set of natural numbers.
For $n\in\mb{N}$ and an indeterminate $x$, define
the rising and falling factorials (cf.~\citu{comtet}{\S1.4})
by $(x)_0=\ang{x}_0=0$ and
\pnq{(x)_n&=x(x+1)\cdots(x+n-1),\\
\ang{x}_n&=x(x-1)\cdots(x-n+1).}
Denote by $[T^m]\phi(T)$ the coefficient (Wilf~\citu{wilf}{\S1.2})
of $T^m$ in the formal power series $\phi(T)$. The following
algebraic relation with symmetric functions as connection
coefficients is crucial in our computations.
\begin{lemm}[Chen and Chu~\cito{chen+w}] \label{sigma}
Let $m$ be a nonnegative integer and $\{x,y\}$
two indeterminates. Then the following algebraic equation holds
\[x^{2m}=\sum_{{\ell}=0}^m(-1)^{\ell}\ang{y+x}_{\ell}\ang{y-x}_{\ell}\sig_{m,\ell}(y),\]
where $\sig_{m,\ell}(y)$ is a complete symmetric function
\pnq{\sig_{m,\ell}(y)&=[T^{m-\ell}]\prod_{\jj=0}^{\ell}\frac1{1-T(y-\jj)^2}\\
&=\sum\sbs{0\le k_1\le k_2\le\cdots\le k_{m-\ell}\le\ell}\prod_{\ii=1}^{m-\ell}(y-k_\ii)^2,}
and admits the explicit expression:
\[\sig_{m,\ell}(y)=\frac{2(-1)^{\ell}}{\ang{2y}_{1+{2\ell}}}
\sum_{i=0}^{\ell}\binm{2y}{i}\binm{{2\ell}-2y}{{\ell}-i}(y-i)^{1+2m}.\]
\end{lemm}

There exist numerous binomial identities in the literature
(cf.~\cite{kn:bai+w,kn:chu02a,kn:elsner,kn:gut+hm,kn:hilton,kn:shapiro,kn:slavik})
and summation formulae for binomial moments
(cf.~\cite{kn:chen+w,kn:chu17a,kn:chu19e,kn:chu20b,kn:miana}).
By combining the telescoping method (cf.~\citu{graham}{\S5.8})
with the above algebraic relation,
we shall examine four classes of binomial moments
\[\sum_{k=1}^n(\pm1)^{k-1}\binm{2n}{n-k}k^m
\quad\text{and}\quad
\sum_{k=1}^n(\pm1)^{k-1}\binq{2n}{n-k}k^m,\]
where $\binm{2n}{n-k}$ is the usual binomial coefficient and
\[\binq{n}{k}=\frac{(\frac12)_n}{(\frac12)_k(\frac12)_{n-k}}
\quad\text{for}\quad n,k\in\mb{N}.\]
Several explicit formulae will be established in the next four sections.
Finally, an unresolved problem will be proposed at the end of the paper.

\section{Positive Sums
\boxed{A_m(n):=\sum_{k=1}^n\binm{2n}{n-k}k^m}}
As a warm-up, we begin with three initial cases.
First, it is almost trivial to evaluate
\[A_0(n)=\frac12\sum_{k=-n}^n\binm{2n}{n-k}-\binm{2n}{n}=2^{2n-1}-\frac12\binm{2n}{n}.\]

Next for $m=1$, by applying \boxed{k=\frac{n+k}{2}-\frac{n-k}{2}}
and telescoping, we derive the closed formula
\pnq{A_1(n)
&=\sum_{k=1}^n\binm{2n}{n-k}\bigg\{\frac{n+k}{2}-\frac{n-k}{2}\bigg\}\\
&=n\sum_{k=1}^n\bigg\{\binm{2n-1}{n-k}-\binm{2n-1}{n-k-1}\bigg\}\\
&=n\binm{2n-1}{n-1}=\frac n2\binm{2n}{n}.}

Then for $m=2$, by writing \boxed{k^2=n^2-(n^2-k^2)},
we deduce the following simpler expression
\pnq{A_2(n)&=n^2A_0(n)-\ang{2n}_2A_0(n-1)\\
&=n^2\bigg\{2^{2n-1}-\frac12\binm{2n}{n}\bigg\}
-2n(2n-1)\bigg\{2^{2n-3}-\frac12\binm{2n-2}{n-1}\bigg\}\\
&=2^{2n-2}n.}

In order to handle with the general case, we let $x\to k$
and $y\to n$ in Lemma~\ref{sigma}:
\pq{\label{sig=g}
k^{2m}=\sum_{{\ell}=0}^m(-1)^{\ell}\ang{n+k}_{\ell}\ang{n-k}_{\ell}\sig_{m,\ell}(n).}
In particular, $\sig_{m,\ell}(n)$ can be expressed explicitly as
\pq{\label{sig+g}
\sig_{m,\ell}(n)=\frac{2(-1)^{\ell}}{\ang{2n}_{1+{2\ell}}}
\sum_{i=0}^{\ell}\binm{2n}{i}\binm{{2\ell}-2n}{{\ell}-i}(n-i)^{1+2m},}
which is valid for $\color{blue}n>m$, avoiding zero factors
in the denominator.

\subsection{Moments of even powers} \
By substitution, we can manipulate the double sum
\pnq{A_{2m}(n)
&=\sum_{k=1}^n\binm{2n}{n-k}k^{2m}=\frac12\sum_{k=-n}^n\binm{2n}{n-k}k^{2m}\\
&=\frac12\sum_{k=-n}^n\binm{2n}{n-k}\sum_{{\ell}=0}^m(-1)^{\ell}\ang{n+k}_{\ell}\ang{n-k}_{\ell}\sig_{m,\ell}(n)\\
&=\frac12\sum_{{\ell}=0}^m(-1)^{\ell}\ang{2n}_{2\ell}\sig_{m,\ell}(n)\sum_{k=-n}^n\binm{2n-{2\ell}}{n-{\ell}-k}\\
&=\sum_{{\ell}=0}^m(-1)^{\ell}2^{2n-{2\ell}-1}\ang{2n}_{2\ell}\sig_{m,\ell}(n).}

Therefore, we have established the formula as in the theorem below.
\begin{thm}[$m,n\in\mb{N}$]\label{thm=Ag}
\[A_{2m}(n)=\sum_{{\ell}=0}^m(-1)^{\ell}2^{2n-{2\ell}-1}\ang{2n}_{2\ell}\sig_{m,\ell}(n).\]
\end{thm}

Applying the explicit expression in \eqref{sig+g} for $\sig_{m,\ell}(n)$,
we can derive, from this theorem, several concrete summation formulae
for small values of $m$. The first five of them are highlighted below.
\begin{corl}
\pnq{A_2(n)&=2^{2 n-2} n,\\
A_4(n)&=2^{2 n-3} n (3 n-1),\\
A_6(n)&=2^{2 n-4} n (15 n^2-15 n+4),\\
A_8(n)&=2^{2 n-5} n (105 n^3-210 n^2+147 n-34),\\
A_{10}(n)&=2^{2 n-6} n (945 n^4-3150 n^3+4095 n^2-2370 n+496).}
\end{corl}

\subsection{Moments of odd powers} \
Analogously, we can deal with the case of odd powers
\pnq{A_{1+2m}(n)
&=\sum_{k=1}^n\binm{2n}{n-k}k^{1+2m}=\sum_{k=1}^nk\binm{2n}{n-k}k^{2m}\\
&=\sum_{k=1}^nk\binm{2n}{n-k}\sum_{{\ell}=0}^m(-1)^{\ell}\ang{n+k}_{\ell}\ang{n-k}_{\ell}\sig_{m,\ell}(n)\\
&=\sum_{{\ell}=0}^m(-1)^{\ell}\ang{2n}_{2\ell}\sig_{m,\ell}(n)\sum_{k=1}^nk\binm{2n-{2\ell}}{n-{\ell}-k}.}
Evaluating the inner sum with respect to $k$
\pnq{\sum_{k=1}^nk\binm{2n-{2\ell}}{n-{\ell}-k}
&=\sum_{k=1}^n\binm{2n-{2\ell}}{n-{\ell}-k}
\bigg\{\frac{n-{\ell}+k}2-\frac{n-{\ell}-k}2\bigg\}\\
=&(n-{\ell})\sum_{k=1}^n\bigg\{\binm{2n-{2\ell}-1}{n-{\ell}-k}-
\binm{2n-{2\ell}-1}{n-{\ell}-k-1}\bigg\}\\
=&(n-{\ell})\binm{2n-{2\ell}-1}{n-{\ell}-1}
=\frac{n-{\ell}}2\binm{2n-{2\ell}}{n-{\ell}},}
and then making substitution, we derive the followng formula
after slight simplifications.
\begin{thm}[$m,n\in\mb{N}$]\label{thm=Ah}
\[A_{1+2m}(n)=\frac12\binm{2n}{n}
\sum_{{\ell}=0}^m(-1)^{\ell}\ang{n}_{\ell}\ang{n}_{\ell+1}\sig_{m,\ell}(n).\]
\end{thm}

Replacing $\sig_{m,\ell}(n)$ by \eqref{sig+g} in this theorem,
we establish the following simplified formulae.
\begin{corl}
\pnq{A_1(n)&=\frac{n}{2}\binm{2n}{n},\\
A_3(n)&=\frac{n^2}{2}\binm{2n}{n},\\
A_5(n)&=\frac{n^2}{2}\binm{2n}{n}(2 n-1),\\
A_7(n)&=\frac{n^2}{2}\binm{2n}{n}(6 n^2-8 n+3),\\
A_9(n)&=\frac{n^2}{2}\binm{2n}{n}(24 n^3-60 n^2+54 n-17).}
\end{corl}

\section{Alternating Sums
\boxed{B_m(n):=\sum_{k=1}^n(-1)^{k-1}\binm{2n}{n-k}k^m}}
First, it is easy to compute
\[B_0(n)=\sum_{k=1}^n(-1)^{k-1}\bigg\{\binm{2n-1}{n-k}+\binm{2n-1}{n-k-1}\bigg\}\\
=\frac12\binm{2n}{n}.\]

Next for $m=1$, we can evaluate the corresponding sum in the same manner as $A_1(n)$
\pnq{B_1(n)&=n\sum_{k=1}^n(-1)^{k-1}\bigg\{\binm{2n-1}{n-k}-\binm{2n-1}{n-k-1}\bigg\}\\
&=n\sum_{k=1}^n(-1)^{k-1}\bigg\{\binm{2n-2}{n-k}+\binm{2n-2}{n-k-2}\bigg\}\\
&=n\bigg\{\binm{2n-2}{n-1}-\binm{2n-2}{n-2}\bigg\}\\
&=\frac{n}{2(2n-1)}\binm{2n}{n}=\binm{2n-2}{n-1}.}

Then for $m=2$, we have by writing \boxed{k^2=n^2-(n^2-k^2)}
\[B_2(n)=n^2B_0(n)-\ang{2n}_2B_0(n-1)=\chi(n=1).\]

\subsection{Moments of even powers} \
By making use of \eqref{sig=g}, we can proceed with
\pnq{B_{2m}(n)
&=\sum_{k=1}^n(-1)^{k-1}\binm{2n}{n-k}k^{2m}=\frac12\sum_{k=-n}^n(-1)^{k-1}\binm{2n}{n-k}k^{2m}\\
&=\frac12\sum_{k=-n}^n(-1)^{k-1}\binm{2n}{n-k}\sum_{{\ell}=0}^m(-1)^{\ell}\ang{n+k}_{\ell}\ang{n-k}_{\ell}\sig_{m,\ell}(n)\\
&=\frac12\sum_{{\ell}=0}^m(-1)^{\ell}\ang{2n}_{2\ell}\sig_{m,\ell}(n)\sum_{k=-n}^n(-1)^{k-1}\binm{2n-{2\ell}}{n-{\ell}-k}.}
Evaluating the rightmost binomial sum
\[\color{blue}\sum_{k=-n}^n(-1)^{k}\binm{2n-{2\ell}}{n-{\ell}-k}
=(-1)^{n-\ell}\binm{\ell}{2n-\ell}\]
and then simplifying the resulting expression,
we establish the following formula.
\begin{thm}[$m,n\in\mb{N}$]\label{thm=Bg}
\[B_{2m}(n)=\frac{(-1)^{n-1}}2\sum_{{\ell}=0}^m\ang{2n}_{2\ell}
\binm{\ell}{2n-\ell}\sig_{m,\ell}(n).\]
\end{thm}
When $n>m$, the binomial coefficient \boxed{\binm{\ell}{2n-\ell}=0}
since $0\le\ell\le m$. From this, we arrive at the simplest value below.
\begin{corl}
\[B_{2m}(n)=0 \quad\text{for}\quad n>m.\]
\end{corl}

\subsection{Moments of odd powers} \
Now, we examine the case of odd powers:
\pnq{B_{1+2m}(n)
&=\sum_{k=1}^n(-1)^{k-1}\binm{2n}{n-k}k^{1+2m}=\sum_{k=1}^n(-1)^{k-1}k\binm{2n}{n-k}k^{2m}\\
&=\sum_{k=1}^n(-1)^{k-1}k\binm{2n}{n-k}\sum_{{\ell}=0}^m(-1)^{\ell}\ang{n+k}_{\ell}\ang{n-k}_{\ell}\sig_{m,\ell}(n)\\
&=\sum_{{\ell}=0}^m(-1)^{\ell}\ang{2n}_{2\ell}\sig_{m,\ell}(n)\sum_{k=1}^n(-1)^{k-1}\binm{2n-{2\ell}}{n-{\ell}-k}k.}
The rightmost sum can be evaluated in closed form
\pnq{\sum_{k=1}^n(-1)^{k-1}&k\binm{2n-{2\ell}}{n-{\ell}-k}
=\sum_{k=1}^n(-1)^{k-1}\bigg\{\frac{n-{\ell}+k}2-\frac{n-{\ell}-k}2\bigg\}\binm{2n-{2\ell}}{n-{\ell}-k}\\
=&(n-{\ell})\sum_{k=1}^n(-1)^{k-1}\bigg\{\binm{2n-{2\ell}-1}{n-{\ell}-k}-
\binm{2n-{2\ell}-1}{n-{\ell}-k-1}\bigg\}\\
=&(n-{\ell})\sum_{k=1}^n(-1)^{k-1}\bigg\{\binm{2n-{2\ell}-2}{n-{\ell}-k}-
\binm{2n-{2\ell}-2}{n-{\ell}-k-2}\bigg\}\\
=&(n-{\ell})\bigg\{\binm{2n-2\ell-2}{n-\ell-1}-\binm{2n-2\ell-2}{n-\ell-2}\bigg\}
=\binm{2n-2\ell-2}{n-\ell-1}.}
By substitution, we arrive at the following summation formula.
\begin{thm}[$m,n\in\mb{N}$]\label{thm=Bh}
\[B_{1+2m}(n)=\sum_{{\ell}=0}^m(-1)^{\ell}\ang{2n}_{2\ell}
\binm{2n-2\ell-2}{n-\ell-1}\sig_{m,\ell}(n).\]
\end{thm}

Recalling further \eqref{sig+g} for $\sig_{m,\ell}(n)$,
we have the five initial identities as below.
\begin{corl}[$n\in\mb{N}$]
\pnq{B_1(n)&=\binm{2n}{n}\frac{n}{2 (2 n-1)},\\
B_3(n)&=\binm{2n}{n}\frac{-n^2}{2 (2 n-1)(2 n-3)},\\
B_5(n)&=\binm{2n}{n}\frac{n^2 (4 n-1)}{2(2 n-1) (2 n-3) (2 n-5)},\\
B_7(n)&=\binm{2n}{n}\frac{-n^2 \big(34 n^2-24 n+5\big) }{2 (2 n-1) (2 n-3) (2 n-5) (2 n-7)},\\
B_9(n)&=\binm{2n}{n}\frac{n^2 \big(496 n^3-672 n^2+344 n-63\big)}{2 (2 n-1) (2 n-3) (2 n-5) (2 n-7) (2 n-9)}.}
\end{corl}

\section{Alternating Sums
\boxed{C_m(n):=\sum_{k=1}^n(-1)^{k-1}\binq{2n}{n-k}k^m}}
According to the recurrence relation
\[\binq{2n}{n-k}=\frac{2n-\frac12}{2n-1}
\bigg\{\binq{2n-1}{n-k}+\binq{2n-1}{n-k-1}\bigg\},\]
we can evaluate the initial sum in closed form by telescoping
\pnq{C_0(n)&=\frac{2n-\frac12}{2n-1}\sum_{k=1}^n
(-1)^{k-1}\bigg\{\binq{2n-1}{n-k}+\binq{2n-1}{n-k-1}\bigg\}\\
&=\frac{2n-\frac12}{2n-1}
\bigg\{\binq{2n-1}{n-1}-(-1)^n\binq{2n-1}{-1}\bigg\}\\
&=\frac{1}2\binq{2n}{n}+\frac{(-1)^n}{4n-2}.}

Next for $m=1$, applying \boxed{k=\frac{n+k-\tfrac12}{2}-\frac{n-k-\tfrac12}{2}},
we can proceed with
\pnq{C_1(n)&=\sum_{k=1}^n(-1)^{k-1}\binq{2n}{n-k}k\\
&=\sum_{k=1}^n(-1)^{k-1}
\bigg\{\frac{n+k-\tfrac12}{2}\binq{2n}{n-k}
-\frac{n-k-\tfrac12}{2}\binq{2n}{n-k}\bigg\}\\
&=\Big(n-\frac14\Big)\sum_{k=1}^n(-1)^{k-1}
\bigg\{\binq{2n-1}{n-k}-\binq{2n-1}{n-k-1}\bigg\}\\
&=\Big(n-\frac14\Big)\frac{4n-3}{4n-4}
\sum_{k=1}^n(-1)^{k-1}\bigg\{\binq{2n-2}{n-k}-\binq{2n-2}{n-k-2}\bigg\}\\
&=\Big(n-\frac14\Big)\frac{4n-3}{4n-4}
\bigg\{\binq{2n-2}{n-1}-\binq{2n-2}{n-2}
-(-1)^n\binq{2n-2}{-1}+(-1)^n\binq{2n-2}{-2}\bigg\}\\
&=\frac{2n-1}{8n-8}\binq{2n}{n}+(-1)^n\frac{2n+1}{8n-8}.}

Then for $m=2$, we have by writing \boxed{k^2=\Big(n-\frac12\Big)^2-\bigg\{\Big(n-\frac12\Big)^2-k^2\bigg\}}
\pnq{C_2(n)&=\Big(n-\frac12\Big)^2C_0(n)
-\Big(2n-\frac12\Big)\Big(2n-\frac32\Big)C_0(n-1)\\
&=\Big(n-\frac12\Big)^2\bigg\{\frac{1}2\binq{2n}{n}+\frac{(-1)^n}{4n-2}\bigg\}\\
&-\Big(2n-\frac12\Big)\Big(2n-\frac32\Big)
\bigg\{\frac{1}2\binq{2n-2}{n-1}+\frac{(-1)^n}{4n-6}\bigg\}\\
&=(-1)^n\frac{n(n+1)}{2(2n-3)}.}

In order to treat the general case, we rewrite the equations in Lemma~\ref{sigma}
by making replacements $x$ $\to$ $k$ and $y$ $\to$ $n-\frac12$:
\pq{\label{sig=h}
k^{2m}=\sum_{{\ell}=0}^m(-1)^{\ell}\ang{n-\tfrac12+k}_{\ell}
\ang{n-\tfrac12-k}_{\ell}\sig_{m,\ell}\big(n-\tfrac12\big).}
In particular, $\sig_{m,\ell}\big(n-\tfrac12\big)$ can be expressed explicitly as
\pq{\label{sig+h}
\sig_{m,\ell}\big(n-\tfrac12\big)=\frac{2(-1)^{\ell}}{\ang{2n-1}_{1+{2\ell}}}
\sum_{i=0}^{\ell}\binm{2n-1}{i}\binm{{2\ell}-2n+1}{{\ell}-i}\big(n-i-\tfrac12\big)^{1+2m},}
which is valid for $\color{blue}n>m$, avoiding zero factors
in the denominator.

\subsection{Moments of even powers} \
By means of linear relation \eqref{sig=h},
we reformulate the double sum
\pnq{C_{2m}(n)&=\sum_{k=1}^n(-1)^{k-1}\binq{2n}{n-k}k^{2m}\\
&=\sum_{k=1}^n(-1)^{k-1}\binq{2n}{n-k}\sum_{{\ell}=0}^m(-1)^{\ell}\ang{n-\tfrac12+k}_{\ell}\ang{n-\tfrac12-k}_{\ell}\sig_{m,\ell}\big(n-\tfrac12\big)\\
&=\sum_{{\ell}=0}^m(-1)^{\ell}\ang{2n-\tfrac12}_{2\ell}\sig_{m,\ell}\big(n-\tfrac12\big)\sum_{k=1}^n(-1)^{k-1}\binq{2n-{2\ell}}{n-{\ell}-k}.}
The rightmost sum with respect to $k$ can be evaluated as
\pnq{&\sum_{k=1}^n(-1)^{k-1}\binq{2n-{2\ell}}{n-{\ell}-k}\\
=&\frac{2n-{2\ell}-\frac12}{2n-{2\ell}-1}\sum_{k=1}^n(-1)^{k-1}
\bigg\{\binq{2n-{2\ell}-1}{n-{\ell}-k}+\binq{2n-{2\ell}-1}{n-{\ell}-k-1}\bigg\}\\
=&\frac{2n-{2\ell}-\frac12}{2n-{2\ell}-1}
\bigg\{\binq{2n-{2\ell}-1}{n-{\ell}-1}-(-1)^{n}\binq{2n-{2\ell}-1}{-{\ell}-1}\bigg\}\\
=&\frac12\bigg\{\binq{2n-{2\ell}}{n-{\ell}}-(-1)^{n}\frac{4n-{2\ell}+1}{2n-{2\ell}-1}\binq{2n-{2\ell}}{-{\ell}-1}\bigg\}.}
By substitution, we deduce the following expression
\[C_{2m}(n)=\frac12
\sum_{{\ell}=0}^m(-1)^{\ell}\ang{2n-\tfrac12}_{2\ell}
\bigg\{\binq{2n-{2\ell}}{n-{\ell}}-(-1)^{n}\frac{4n-{2\ell}+1}{2n-{2\ell}-1}
\binq{2n-{2\ell}}{-{\ell}-1}\bigg\}\sig_{m,\ell}\big(n-\tfrac12\big).\]
Furthermore, we have the following unexpected identity
\[\color{blue}\sum_{{\ell}=0}^m(-1)^{\ell}\ang{2n-\tfrac12}_{2\ell}
\binq{2n-{2\ell}}{n-{\ell}}\sig_{m,\ell}\big(n-\tfrac12\big)=0
\quad\text{for}\quad m,n\in\mb{N}.\]
In fact, we can rewrite the sum $\Lam(m,n)$ on the left according to Lemma~\ref{sigma}
\pnq{\Lam(m,n)&=\sum_{{\ell}=0}^m(-1)^{\ell}\ang{2n-\tfrac12}_{2\ell}
\binq{2n-{2\ell}}{n-{\ell}}\sig_{m,\ell}\big(n-\tfrac12\big)\\
&=\big[T^{m}\big]\sum_{{\ell}=0}^m(-1)^{\ell}
\frac{(\tfrac12)_{2n}~T^\ell}{(\frac12)^2_{n-\ell}
\prod_{\jj=0}^{\ell}\{1-T(n-\jj-\frac12)^2\}}.}  

For the sequence $\lam_\ell$ defined below
\[\lam_{\ell}:=\frac{(\tfrac12)_{2n}~T^\ell}{(\frac12)^2_{n-\ell}
\prod_{\jj=0}^{\ell-1}\{1-T(n-\jj-\frac12)^2\}},\]
we can compute the sum of two consecutive terms 
\[\lam_{\ell}+\lam_{\ell+1}
=\frac{(\tfrac12)_{2n}~T^\ell}{(\frac12)^2_{n-\ell}
\prod_{\jj=0}^{\ell}\{1-T(n-\jj-\frac12)^2\}}.\]
Thus, we can further evaluate for $m,n\in\mb{N}$
\pnq{\Lam(m,n)&=\big[T^{m}\big]
\bigg\{\frac{(\tfrac12)_{2n}}{(\frac12)^2_{n}\big(1-T(n-\frac12)^2\big)}
+\sum_{{\ell}=1}^m(-1)^{\ell}\big(\lam_{\ell}+\lam_{\ell+1}\big)\bigg\}\\
&=\big[T^{m}\big]
\bigg\{\frac{(\tfrac12)_{2n}}{(\frac12)^2_{n}\big(1-T(n-\frac12)^2\big)}
-\lam_1+(-1)^m\lam_{m+1}\bigg\}\\
&=\big[T^{m}\big]\binq{2n}{n}=0.}

Consequently, we establish the following formula.
\begin{thm}[$m,n\in\mb{N}$]\label{thm=Cg}
\pnq{C_{2m}(n)&=\sum_{{\ell}=0}^m(-1)^{n}\frac{4n-2\ell+1}{4n-4\ell-2}
\frac{\ang{2n-\tfrac12}_{2\ell}}{\binq{2n-\ell+1}{\ell+1}}\sig_{m,\ell}\big(n-\tfrac12\big)\\
&=\sum_{\ell=0}^m(-1)^{n}\frac{(\tfrac12)_{\ell}(\tfrac12)_{\ell+1}}{2n-2\ell-1}
\binq{2n}{\ell}\sig_{m,\ell}\big(n-\tfrac12\big).}
\end{thm}

According to \eqref{sig+h}, we can reduce the above formulae further
for small values of $m$ that gives rise to the following five identities.
\begin{corl}
\pnq{C_2(n)&=\frac{(-1)^n n (n+1)}{2 (2 n-3)} \\
C_4(n)&=\frac{(-1)^n n (n+1) }{2 (2 n-3) (2 n-5)}\big\{2 n^3-n^2-5 n+1\big\} \\
C_6(n)&=\frac{(-1)^n n (n+1)}{2 (2 n-3) (2 n-5) (2 n-7)}
 \big\{4 n^6-8 n^5-25 n^4+30 n^3+40 n^2-31 n+5\big\} \\
C_8(n)&=\frac{(-1)^n n (n+1) }{2 (2 n-3) (2 n-5) (2 n-7) (2 n-9)}
\Bigg\{\mult{c}{8 n^9-36 n^8-62 n^7+301 n^6+231 n^5\\
-847 n^4-175 n^3+855 n^2-443 n+63}\Bigg\}\\
C_{10}(n)&\pp{t}{&=\frac{(-1)^n n (n+1) }{2 (2 n-3) (2 n-5) (2 n-7) (2 n-9) (2 n-11)}\\
&\times\Bigg\{\mult{c}{16 n^{12}-128 n^{11}-8 n^{10}+1680 n^9-735 n^8-9348 n^7+4368 n^6\\
+23466 n^5-17070 n^4-19460 n^3+28666 n^2-12077 n+1575}\Bigg\}.}}
\end{corl}

\subsection{Moments of odd powers} \
We can deal with the case of odd powers similarly:
\pnq{C_{1+2m}(n)
&=\sum_{k=1}^n(-1)^{k-1}\binq{2n}{n-k}k^{1+2m}\\
&=\sum_{k=1}^n(-1)^{k-1}k\binq{2n}{n-k}\sum_{{\ell}=0}^m(-1)^{\ell}
\ang{n-\tfrac12+k}_{\ell}\ang{n-\tfrac12-k}_{\ell}\sig_{m,\ell}(n-\tfrac12)\\
&=\sum_{{\ell}=0}^m(-1)^{\ell}\ang{2n-\tfrac12}_{2\ell}\sig_{m,\ell}(n-\tfrac12)
\sum_{k=1}^n(-1)^{k-1}\binq{2n-{2\ell}}{n-{\ell}-k}k.}
Evaluating the above sum with respect to $k$
\pnq{&\sum_{k=1}^n(-1)^{k-1}\binq{2n-{2\ell}}{n-{\ell}-k}k
=\sum_{k=1}^n(-1)^{k-1}\binq{2n-{2\ell}}{n-{\ell}-k}
\bigg\{\frac{n-\ell+k-\tfrac12}2-\frac{n-\ell-k-\tfrac12}2\bigg\}\\
=&\frac{2n-{2\ell}-\tfrac12}{2}\sum_{k=1}^n(-1)^{k-1}\bigg\{\binq{2n-{2\ell}-1}{n-{\ell}-k}-
\binq{2n-{2\ell}-1}{n-{\ell}-k-1}\bigg\}\\
=&\frac{(2n-{2\ell}-\tfrac12)(2n-{2\ell}-\tfrac32)}{2(2n-{2\ell}-2)}
\sum_{k=1}^n(-1)^{k-1}\bigg\{\binq{2n-{2\ell}-2}{n-{\ell}-k}-
\binq{2n-{2\ell}-2}{n-{\ell}-k-2}\bigg\}\\
=&\pp{t}{\frac{(2n-{2\ell}-\tfrac12)(2n-{2\ell}-\tfrac32)}{2(2n-{2\ell}-2)}
\bigg\{&\binq{2n-{2\ell}-2}{n-{\ell}-1}-\binq{2n-{2\ell}-2}{n-{\ell}-2}\\
-(-1)^n&\binq{2n-{2\ell}-2}{-{\ell}-1}+(-1)^n\binq{2n-{2\ell}-2}{-{\ell}-2}\bigg\}}\\
=&\frac{2n-2\ell-1}{8(n-\ell-1)}\binq{2n-{2\ell}}{n-{\ell}}+(-1)^n\frac{(1+2n)(1+2\ell)}{8(n-{\ell}-1)}
\binq{2n-{2\ell}}{-{\ell}},}
we obtain the following explicit formula.
\begin{thm}[$n>m+1$ with $m,n\in\mb{N}$] \label{thm=Ch}
\pnq{C_{1+2m}(n)&=\frac18\sum_{{\ell}=0}^m(-1)^{\ell}\ang{2n-\tfrac12}_{2\ell}
\bigg\{\frac{2n-2\ell-1}{n-{\ell}-1}\binq{2n-{2\ell}}{n-{\ell}}+(-1)^n
\frac{(1+2n)(1+2\ell)}{n-{\ell}-1}
\binq{2n-{2\ell}}{-{\ell}}\bigg\}\sig_{m,\ell}(n)\\
&=\frac18\sum_{{\ell}=0}^m(-1)^{\ell}\ang{2n-\tfrac12}_{2\ell}
\bigg\{\frac{2n-2\ell-1}{n-{\ell}-1}\binq{2n-{2\ell}}{n-{\ell}}
+(-1)^{n+\ell}\frac{(1+2n)(1+2\ell)}{n-{\ell}-1}
\binq{2n-\ell}{\ell}^{-1}\bigg\}\sig_{m,\ell}(n).}
\end{thm}

By substituting \eqref{sig+h} into the above sums
and then simplifying the resulting expressions,
we find the following five summation formulae.
\begin{corl}
\pnq{\boxed{n>1}\quad
C_{1}(n)&=(-1)^n\frac{(2 n+1)}{8 (n-1)}+\binq{2n}{n}\frac{(2 n-1)}{8 (n-1)},\\[2mm]
\boxed{n>2}\quad
C_{3}(n)&=(-1)^n\frac{(2 n+1) (4 n^3-6 n+1)}{32 (n-1) (n-2)}
-\binq{2n}{n}\frac{(2 n-1)^2}{32 (n-1) (n-2)},\\[2mm]
\boxed{n>3}\quad
C_{5}(n)&=(-1)^n\tfrac{(2 n+1) (8 n^6-8 n^5-40 n^4+20 n^3+40 n^2-22 n+3)}{64 (n-1) (n-2)(n-3)}
+\binq{2n}{n}\tfrac{(2 n-1)^2 (4 n-3)}{64 (n-1) (n-2)(n-3)},\\[3mm]
\boxed{n>4}\quad
C_{7}(n)&=(-1)^n\tfrac{(2 n+1) (32 n^9-96 n^8-272 n^7+616 n^6+840 n^5-1288 n^4-532 n^3+1068 n^2-422 n+51)}{256 (n-1) (n-2)(n-3)(n-4)}\\
&-\binq{2n}{n}\frac{(2 n-1)^2 (68 n^2-116 n+51)}{256 (n-1) (n-2)(n-3)(n-4)},\\[3mm]
\boxed{n>5}\quad
C_{9}(n)&=(-1)^n\tfrac{(2 n+1)
\resizebox{9cm}{!}{$\bigg\{\mult{c}{32 n^{12}-192 n^{11}-224 n^{10}+2208 n^9+864 n^8-9744 n^7-840 n^6\\
+18792 n^5-7224 n^4-12532 n^3+12576 n^2-4178 n+465}\bigg\}$}}{256 (n-1) (n-2)(n-3)(n-4)(n-5)}\\
&+\binq{2n}{n}\frac{(2 n-1)^2 (496 n^3-1416 n^2+1388 n-465)}{256 (n-1) (n-2)(n-3)(n-4)(n-5)}.}
\end{corl}

\section{Positive Sums
\boxed{D_m(n):=\sum_{k=1}^n\binq{2n}{n-k}k^m}}

First for $m=1$, by writing
\[k=\frac{n+k-\tfrac12}{2}-\frac{n-k-\tfrac12}{2},\]
we can calculate the corresponding sum as follows:
\pnq{D_1(n)&=\sum_{k=1}^nk\binq{2n}{n-k}
=\sum_{k=1}^n\bigg\{\frac{n+k-\tfrac12}2\binq{2n}{n-k}
-\frac{n-k-\tfrac12}2\binq{2n}{n-k}\bigg\}\\
&=\Big(n-\frac14\Big)\sum_{k=1}^n
\bigg\{\binq{2n-1}{n-k}-\binq{2n-1}{n-k-1}\bigg\}\\
&=\Big(n-\frac14\Big)\frac{4n-3}{4n-4}
\sum_{k=1}^n\bigg\{\binq{2n-2}{n-k}-\binq{2n-2}{n-k-2}\bigg\}\\
&=\Big(n-\frac14\Big)\frac{4n-3}{4n-4}
\bigg\{\binq{2n-2}{n-1}+\binq{2n-2}{n-2}
-\binq{2n-2}{-1}-\binq{2n-2}{-2}\bigg\}\\
&=\frac14+\binq{2n}{n}\frac{2n-1}{4}.}

Then we turn to examine, in general, the corresponding sums
to the case of odd powers:
\pnq{D_{1+2m}(n)
&=\sum_{k=1}^n\binq{2n}{n-k}k^{1+2m}=\sum_{k=1}^nk\binq{2n}{n-k}k^{2m}\\
&=\sum_{k=1}^nk\binq{2n}{n-k}\sum_{{\ell}=0}^m(-1)^{\ell}\ang{n-\tfrac12+k}_{\ell}
\ang{n-\tfrac12-k}_{\ell}\sig_{m,\ell}\big(n-\tfrac12\big)\\
&=\sum_{{\ell}=0}^m(-1)^{\ell}\ang{2n-\tfrac12}_{2\ell}
\sig_{m,\ell}\big(n-\tfrac12\big)\sum_{k=1}^nk\binq{2n-{2\ell}}{n-{\ell}-k}.}
The sum with respect to $k$ can be evaluated as follows:
\pnq{&\sum_{k=1}^nk\binq{2n-{2\ell}}{n-{\ell}-k}
=\sum_{k=1}^n\binq{2n-{2\ell}}{n-{\ell}-k}
\bigg\{\frac{n-{\ell}+k-\tfrac12}2-\frac{n-{\ell}-k-\tfrac12}2\bigg\}\\
&~=\frac{4n-4\ell-1}{4}\sum_{k=1}^n\bigg\{\binq{2n-{2\ell}-1}{n-{\ell}-k}-
\binq{2n-{2\ell}-1}{n-{\ell}-k-1}\bigg\}\\
&~=\frac{(4n-4\ell-1)(4n-4\ell-3)}{16(n-\ell-1)}
\sum_{k=1}^n\bigg\{\binq{2n-{2\ell}-2}{n-{\ell}-k}
-\binq{2n-{2\ell}-2}{n-{\ell}-k-2}\bigg\}\\
&~=\frac{(4n-4\ell-1)(4n-4\ell-3)}{16(n-\ell-1)}
\pp{t}{\bigg\{&\binq{2n-{2\ell}-2}{n-{\ell}-1}+\binq{2n-{2\ell}-2}{n-{\ell}-2}\\
-&\binq{2n-{2\ell}-2}{-{\ell}-1}-\binq{2n-{2\ell}-2}{-{\ell}-2}\bigg\}}\\
&~=\frac{n-{\ell}-\tfrac12}{2}\binq{2n-{2\ell}}{n-{\ell}}+\frac{{\ell}+\tfrac12}{2}
\binq{2n-{2\ell}}{-{\ell}}.}
By substitution, we obtain the following summation formula.
\begin{thm}[$m,n\in\mb{N}$] \label{thm=Dh}
\pnq{D_{1+2m}(n)
&=\sum_{{\ell}=0}^m\ang{2n-\tfrac12}_{{2\ell}}\bigg\{
\frac{2n-2\ell-1}{4(-1)^{\ell}}\binq{2n-{2\ell}}{n-{\ell}}
+\frac{1+2\ell}{4}\binq{2n-\ell}{\ell}^{-1}\bigg\}\sig_{m,\ell}\big(n-\tfrac12\big).}
\end{thm}

Taking into account \eqref{sig+h}, we can further show, from this theorem,
the five identities below.
\begin{corl} \label{Dh}
\pnq{D_1(n)&=\binq{2n}{n}\frac{(2 n-1) }{4}
+\frac{1}{4},\\
D_3(n)&=\binq{2n}{n}\frac{(2 n-1)^2}{8}
+\frac{1}{8} \big(2 n^2+4 n-1\big),\\
D_5(n)&=\binq{2n}{n}\frac{(2 n-1)^2}{4}(n-1)
+\frac{1}{4} \big(n^4+4 n^3+3 n^2-5 n+1\big),\\
D_7(n)&=\binq{2n}{n}\frac{(2n-1)^2}{16}\big(12 n^2-28 n+17\big)\\
&+\frac{1}{16} \big(4 n^6+24 n^5+42 n^4-28 n^3-108 n^2+96 n-17\big),\\
D_{9}(n)&=\binq{2n}{n}\frac{(2 n-1)^2}{4}\big(12 n^3-48 n^2+66 n-31\big)\\
&+\frac{1}{4} \big(n^8+8 n^7+22 n^6+2 n^5-98 n^4-52 n^3+283 n^2-190 n+31\big).}
\end{corl}

\section{Concluding Comments and Problems}
By making use of the algebraic relation stated in Lemma~\ref{sigma},
we have succeeded in showing explicit summation formulae for four classes
of binomial moments. However, for the sums $D_{2m}(n)$ corresponding to
the cases of even powers, the authors fail to determine related explicit
formulae like those in Theorem~\ref{thm=Dh} and Corollary~\ref{Dh}.
Any attempt to resolve this problem is enthusiastically encouraged.


\subsection*{Declarations Competing interests}
The authors declare no competing interests.

\subsection*{Author contributions}
All authors reviewed and approved the manuscript.

\subsection*{Data Availability Statement}
No datasets were generated or analyzed during the current study.
 


\end{document}